\documentclass[11pt]{amsart}
\usepackage[english]{babel}
\usepackage{inputenc}
\usepackage{amsmath,latexsym,amsthm}
\usepackage{graphicx}
\usepackage{amssymb}
\begin{document}
%\begin{flushright}
%\begin{tabular}{l}
%{\sf Uzbek  Mathematical}\\
%{\sf Journal, 2011,\No 2, pp.\pageref{Akhmed1}-\pageref{Akhmed2}}\\
%\end{tabular}
%\end{flushright}

\title [OPTIMAL QUADRATURE FORMULAS]
{Optimal quadrature formulas for the Cauchy type singular integral in the Sobolev space $L_2^{(2)}(0,1)$}

\author {D.M. Akhmedov}

\address{D.M. Akhmedov\\ Institute of mathematics and information technologies,
Tashkent, Uzbekistan.}
\email{hayotov@mail.ru, axmedovdilshod@mail.ru}

\begin{abstract}
In the present paper in $L_2^{(2)}(0,1)$ S.L.Sobolev space the optimal quadrature formula is constructed for
approximate calculation of Cauchy type singular integral.

\textbf{MSC:} 65D32.

\emph{Keywords:} optimal quadrature formulas; the error
functional; the extremal function; Hilbert space; optimal
coefficients.
\end{abstract}

\maketitle

\textbf{1. Introduction}

   Many problems of sciences and technics is naturally reduced to singular integral equations. Moreover (see. \cite{Lif95})  plane problems are reduced to one dimensional singular equations. The theory of one-dimensional singular integ\-ral equations is given in \cite{Gakh77,Muskh68}.  It is known that the solutions of such integral equations are expressed by singular integrals. Therefore approximate cal\-cu\-lation of the singular integrals  with high exactness is actual problem of numerical analysis.

For the singular integral of Cauchy type $\int_0^1\frac{\varphi(x)}{x-t}dx$  we consider the following quadrature formula
$$
\int\limits_0^1\frac{\varphi(x)}{x-t}dx\cong\sum\limits^N_{\beta=0}
C_{\beta}\varphi(x_{\beta}),\eqno(1.1)
$$
where  $0<t<1,$ $C_{\beta}$ are the coefficients, $x_{\beta}$ are the nodes, $N=2,3,4...$, $\varphi$ is a function of the space $L_2^{(m)}(0,1)$.
Here $L_2^{(m)}(0,1)$ is the Sobolev space of functions with a square integrable $m$th generalized derivative.

To the problem of approximate integration of the Cauchy type singular integrals are devoted many works (see, for instance, \cite{Lif95,BelLif85,Gab75,IsrShad91,Shad87} and references therein).

In the work \cite{Shad99} in the space   for the coefficients of the weight optimal quadrature formulas of the form
$$
\int\limits_0^1p(x)\varphi(x)\cong\sum\limits_{\beta=0}^NC[\beta]\varphi(h\beta)\eqno(1.2)
$$
the following system of linear equations is obtained
$$
\sum\limits_{\gamma=0}^NC[\gamma]\cdot\frac{|h\beta-h\gamma|^{2m-1}}{2(2m-1)!}+P_{m-1}[\beta]=f[\beta],\ \ \ \ \ \ [\beta]\in[0,1],\eqno(1.3)
$$
$$
\sum\limits_{\beta=0}^NC[\beta]\cdot[\beta]^{\alpha}=\int\limits_0^1p(x)x^{\alpha}dx, \ \ \ \ \ \ \ \
\alpha=0,1,...,m-1,\eqno(1.4)
$$
where $f[\beta]=\int\limits_0^1p(x)\frac{|x-h\beta|^{2m-1}}{2(2m-1)!}dx,$ $p(x)$ is a weight function,
$[\beta]=h\beta,$ $h=\frac{1}{N},$
$C[\beta]$ are unknown coefficients, $P_{m-1}[\beta]$ is a polynomial of degree  $m-1$.

In particular,  from (1.2) for $p(x)=\frac{1}{x-t},$ $0<t<1$ we get quadrature formula  (1.1) when $x_{\beta}=h\beta=[\beta]$ and
correspondingly from system  (1.3)-(1.4) when $m=2,$ $p(x)=\frac{1}{x-t},$  $0<t<1$ we have the following system of linear equations for the coefficients of optimal  quadrature formula in the sense of Sard in the form (1.1) when $x_{\beta}=h\beta$ in the space  $L_2^{(2)}(0,1)$
$$
\sum\limits_{\gamma=0}^NC([\gamma],t)\cdot\frac{|h\beta-h\gamma|^3}{12}+p_1\cdot[\beta]+p_0=f([\beta],t),\ \ \ \
\ \beta=0,1,2,...,N,\eqno(1.5)
$$
$$
\sum\limits_{\gamma=0}^NC([\gamma],t)=g_0,\eqno(1.6)
$$
$$
\sum\limits_{\gamma=0}^NC([\gamma],t)\cdot [\gamma]=g_1,\eqno(1.7)
$$
where
$$
\begin{array}{ll}
f(h\beta,t)=&\frac{1}{12}\int\limits_0^1\frac{|x-h\beta|^3}{x-t}dx=\frac{1}{12}\bigg(
-\frac{11}{3}(h\beta)^3+(5t+3)(h\beta)^2-\\
&-(2t^2+3t+1,5)(h\beta)+(t^2+\frac{t}{2}+\frac{1}{3})+\\
& +(t-h\beta)^3(-2\ln|h\beta-t|+\ln(t-t^2))\bigg),
\end{array}
\eqno(1.8)
$$
$$
g_0=\int\limits_0^1\frac{dx}{x-t}=\ln\frac{1-t}{t},\eqno(1.9)
$$
$$
g_1=\int\limits_0^1\frac{x}{x-t}dx=1+t\ln\frac{1-t}{t}.\eqno(1.10)
$$
$C([\gamma],t),$  $\gamma=\overline{0,N}$ and $p_1$, $p_0$  are unknowns.

The rest of the paper is organized as follows. In Section 2 we give some auxiliary results and definitions.
Section 3 is devoted to calculation of the optimal coefficients, i.e. we solve system (1.5)-(1.7).

\medskip

\newpage

\textbf{2. Auxiliary results}

\medskip

In the process of solution of system (1.5)-(1.7) the Sobolev method of construction of optimal quadrature formulas (see \cite{Sob74}) is used. In addition, here we use the concept of discrete argument functions and operations on them. The theory of discrete argument functions is given in \cite{Sob74}. For completeness we give some definitions.

\medskip

\textbf{Defnition 2.1.} The function $\varphi[\beta]=\varphi(h\beta)$  is a \emph{function of discrete argument} if it is given on some set of integer values of $\beta$.

\medskip

\textbf{Defnition 2.2.} \emph{The inner product} of two discrete functions $\varphi[\beta]$   and $\psi[\beta]$  is given by
$$
[\varphi,\psi]=\sum\limits_{\beta\in B}\varphi[\beta]\cdot\psi[\beta],\eqno(2.1)
$$
if the series on the right hand side of the last equality converges absolutely.

\medskip

\textbf{Defnition 2.3.} The convolution of two functions $\varphi[\beta]$  and $\psi[\beta]$  is the inner product
$$
\chi[\beta]=\varphi[\beta]\ast\psi[\beta]=[\varphi[\gamma],
\psi[\beta-\gamma]]=\sum\limits_{\gamma=-\infty}^{\infty}\varphi[\gamma]\psi[\beta-\gamma].\eqno(2.2)
$$

Moreover, we need the discrete analogue   $D_2(h\beta)$ of the differential operator  $d^4/dx^4$, which is defined by the following formula (see \cite{Shad85})
$$
D_2(h\beta)=\frac{3!}{h^4}\left\{\begin{array}{ll}
              3!\sqrt{3}q^{|\beta|}, & |\beta|\geq 2,\\
             19-12\sqrt{3},& |\beta|=1, \\
              6\sqrt{3}-8, & \beta=0, \\
            \end{array}\right.\eqno(2.3)
$$
where  $q=\sqrt{3}-2$.

We give some properties of the discrete function $D_2(h\beta)$  from \cite{Sob74}:

$$
D_2(h\beta)\ast(h\beta)^{\alpha}=0,\ \ \ \ 0\leq\alpha<4, \eqno (2.4)
$$
$$
hD_2(h\beta)\ast\frac{|h\beta|^3}{12}=\delta(h\beta), \eqno(2.5)
$$
where  $\delta(h\beta)$ is the discrete delta function and  $\delta(h\beta)=\left\{\begin{array}{cc}
                  0, & \beta\neq 0, \\
                     1, & \beta=0. \\
                   \end{array}  \right.$

\medskip

\textbf{3. Main result}

\medskip

In the present section we solve system (1.5)-(1.7).

 Suppose $C([\beta],t)=0$ when $\beta<0$ and $\beta>N.$
Then, using Definition 2.3, we rewrite system(1.5)-(1.7) in the following convolution form
$$
C([\beta],t)\ast\frac{|h\beta|^3}{12}+p_1\cdot (h\beta)+p_0=f(h\beta,t),\ \
\beta=\overline{0,N},\eqno(3.1)
$$
$$
\sum\limits_{\beta=0}^NC([\beta],t)=g_0,\eqno(3.2)
$$
$$
\sum\limits_{\beta=0}^NC([\beta],t)\cdot(h\beta)=g_1,\eqno(3.3)
$$
where $f(h\beta)$, $g_0$ and $g_1$ are defined by (1.8), (1.9) and (1.10), respectively.

We denote
$$
v(h\beta)=C([\beta],t)\ast\frac{|h\beta|^3}{12},\eqno(3.4)
$$
$$
u(h\beta)=v(h\beta)+p_1\cdot (h\beta)+p_0.\eqno(3.5)
$$
Using the properties (2.4), (2.5) of the operator $D_2(h\beta)$  from (2.3) and (3.5) we have
$$
C([\beta],t)=hD_2(h\beta)\ast u(h\beta).\eqno(3.6)
$$

But for calculation of the convolution  (3.6) we need to define the function  $u(h\beta)$ in all integer values of  $\beta$.
When  $\beta=0,1,2,...,N$ from (3.1) we have  $u(h\beta)=f(h\beta,t)$. Therefore it is sufficient to determine the function
$u(h\beta)$ when $\beta<0$ and $\beta>N$.

Furthere, we define the form of the function   $u(h\beta)$ for $\beta\leq 0$ and $\beta\geq N$. From (3.4) taking into account (3.2), (3.3) we have
$$
v(h\beta)=\left\{\begin{array}{ll}
 -\frac{1}{12}(h\beta)^3g_0 +\frac{1}{4}(h\beta)^2g_1-(h\beta)\cdot p_1^{(0)}-p_0^{(0)},& \beta\leq 0, \\[2mm]
\frac{1}{12}(h\beta)^3g_0 -\frac{1}{4}(h\beta)^2g_1+(h\beta)\cdot p_1^{(0)}+p_0^{(0)}, & \beta\geq N. \\
\end{array}\right.
$$
into account the last equality, from  (3.5) for $u(h\beta)$ we obtain
$$
u(h\beta)=\left\{\begin{array}{ll}
-\frac{1}{12}(h\beta)^3g_0+\frac{1}{4}(h\beta)^2g_1+a_1^-(h\beta)+a_0^-, & \beta\leq 0,\\[2mm]
f(h\beta,t),& 0\leq \beta\leq N, \\[2mm]
\frac{1}{12}(h\beta)^3g_0-\frac{1}{4}(h\beta)^2g_1+a_1^+(h\beta)+a_0^+ ,& \beta\geq N, \\
\end{array}\right.\eqno(3.7)
$$
where $a_1^-$, $a_0^-$, $a_1^+$, $a_0^+$ are unknowns and
$$
\begin{array}{l}
a_1^-=p_1-p_1^{(0)},\ \ a_0^-=p_0-p_0^{(0)},\\
a_1^+=p_1+p_1^{(0)},\ \  a_0^+=p_0+p_0^{(0)}.\\
\end{array}
\eqno (3.8)
$$

Thus, keeping in mind  (2.3), (3.6) and (3.7) we get the following problem

\medskip

\textbf{Problem.} {\it Find the solution of the equation
$$
D_2(h\beta)*u(h\beta)=0 \mbox{ where }  \beta<0, \beta>N\eqno (3.9)
$$
in the form  (3.7)}.

\medskip

If we find unknowns   $a_1^-$, $a_0^-$, $a_1^+$, $a_0^+$, then from  (3.8) we have
$$
\begin{array}{l}
p_1=\frac{1}{2}(a_1^-+a_1^+),\ \ \ \ p_0=\frac{1}{2}(a_0^-+a_0^+),\\[2mm]
p_1^{(0)}=\frac{1}{2}(a_1^+-a_1^-),\ \ p_0^{(0)}=\frac{1}{2}(a_0^+-a_0^-).\\
\end{array}
\eqno (3.10)
$$
Unknowns $a_1^-$, $a_0^-$, $a_1^+$, $a_0^+$ will be found from (3.9), using (2.3) and (3.7).
Then we get the explicit form of the function $u(h\beta)$ and from  (3.6) we find optimal coefficients
$C([\beta],t)$. Moreover, from  (3.10) can be found $p_1,\ p_0$. Thus, the stated problem will be solved completely.

The following holds

\medskip

\textbf{Theorem 1.} \emph{Assume  $ t\neq h\beta$, $\beta=\overline{0,N}$. Then the coefficients of the optimal quadrature formula  (1.1) in the space  $L_2^{(2)}(0,1)$ have the form }
$$
\begin{array}{lcl}
C([0],t)&=&\frac{6}{h^3}\bigg[ \frac{g_0}{12}\left(h^3-3q^N(h^2+h(q+2))\right)+\frac{g_1q^N}{4}(h^2+2h(q+2))+\\
    & & +a_1^-h(q+1)+f(0,t)(3q+2)-f(h,t)(12q+5)-q^N(3f(1,t)\times\\
    & & \times(q+1)+a_1^+h(q+2))+6(q+2)\sum\limits_{\gamma=2}^Nq^{\gamma}f(h\gamma,t) \bigg],
\end{array}\eqno(3.11)
$$
$$
\begin{array}{lcl}
C([\beta],t) &=& \frac{6}{h^3}\bigg[ 6(q+2)\sum\limits_{\gamma=0}^{\beta-2}q^{\beta-\gamma}f(h\gamma,t)-(12q+5)\bigg(f(h(\beta-1),t)+\\
         & & +f(h(\beta+1),t)\bigg)+(6q+4)f(h\beta,t)+6(q+2)\sum\limits_{\gamma=\beta+2}^Nq^{\gamma-\beta}f(h\gamma,t)+\\
         & &+\frac{g_1}{4}\bigg(q^{N-\beta}(2h(q+2)+h^2)-q^{\beta}h^2\bigg)+q^{\beta}\bigg(a_1^-h(q+2)-3f(0,t)\times\\
         & &\times(q+1)\bigg)-q^{N-\beta}(3f(1,t)(q+1)+\frac{g_0}{4}(h(q+2)+h^2)+\\
         & & +a_1^+h(q+2))\bigg], \ \ \ \ \ \ \beta=\overline{1,N-1},
\end{array}\eqno(3.12)
$$
$$
\begin{array}{lcl}
C([N],t)&=& \frac{6}{h^3}\bigg[ -\frac{g_0}{12}(3h(q+1)-h^3)+\frac{g_1}{4}(2h(q+1)-q^Nh^2)+q^N\bigg(a_1^-h(q+2)-\\
     & &-3f(0,t)(q+1)\bigg)-a_1^+h(q+1)+f(1,t)(3q+2)-f(1-h,t)\times\\
     & & \times(12q+5)+6(q+2)\sum\limits_{\gamma=0}^{N-2}q^{N-\gamma}f(h\gamma,t)\bigg].
     \end{array}
\eqno(3.13)
$$
\emph{where}
$$
a_1^-=\frac{\Delta_1}{\Delta}, \eqno (3.14)
$$
$$
a_1^+= \frac{\Delta_2}{\Delta},\eqno (3.15)
$$
$$
a_0^-= f(0,t),\eqno(3.16)
$$
$$
a_0^+=f(1,t)-\frac{1}{12}g_0+\frac{1}{4}g_1-\frac{\Delta_2}{\Delta},\eqno (3.17)
$$
$$
\begin{array}{lcl}
\Delta&=&B_2^2-A_2^2,\\
\Delta_1&=&A_2\bigg[ -F_1-\frac{1}{12}g_0(B_1+3B_2+3B_3)+\frac{1}{4}g_1(B_1+2B_2+B_3-A_3)-\\
        & &-A_1f(0,t)-B_1\bigg(f(1,t)-\frac{1}{12}g_0+\frac{1}{4}g_1\bigg)\bigg]-\\
        & &-B_2\bigg[ -F_2-\frac{1}{12}g_0(A_1+3A_2+3A_3)+\frac{1}{4}g_1(A_1+2A_2+A_3-B_3)-\\
        & &-B_1f(0,t)-A_1\bigg(f(1,t)-\frac{1}{12}g_0+\frac{1}{4}g_1\bigg)\bigg],\\
\Delta_2&=&-A_2\bigg[ -F_2-\frac{1}{12}g_0(A_1+3A_2+3A_3)+\frac{1}{4}g_1(A_1+2A_2+A_3-B_3)-\\
        & &-B_1f(0,t)-A_1\bigg(f(1,t)-\frac{1}{12}g_0+\frac{1}{4}g_1\bigg)\bigg]+\\
        & &+B_2\bigg[ -F_1-\frac{1}{12}g_0(B_1+3B_2+3B_3)+\frac{1}{4}g_1(B_1+2B_2+B_3-A_3)-\\
        & &-A_1f(0,t)-B_1\bigg(f(1,t)-\frac{1}{12}g_0+\frac{1}{4}g_1\bigg)\bigg],\\
F_1&=&6(q+2)\sum\limits_{\gamma=1}^Nq^{\gamma+1}f(h\gamma,t)-(12q+5)f(0,t),\\
F_2&=&6(q+2)\sum\limits_{\gamma=0}^{N-1}q^{N+1-\gamma}f(h\gamma,t)-(12q+5)f(1,t),\\
\end{array}
$$
$$
\begin{array}{lcllcl}
A_1&=&3q+2,\ \    &B_1&=&3q^N(3q+1)\\
A_2&=&h(2q+1),\ \ &B_2&=&hq^N(2q+1)\\
A_3&=&h^2q,\ \    &B_3&=&-h^2q^{N+1},\\
\end{array}
$$
$q=\sqrt{3}-2$.

\medskip

\emph{Proof.} From (3.7) for  $\beta=0$ and $\beta=N$ we get  (3.16) and (3.17), i.e.
$$
\begin{array}{lcl}
a_0^-&=&f(0,t),\\
\end{array}\eqno(3.18)
$$
$$
\begin{array}{c}
a_0^+=f(1,t)-\frac{1}{12}g_0+\frac{1}{4}g_1-a_1^+.\\
\end{array}\eqno(3.19)
$$
From (3.9), using (2.3) and (3.7) for $\beta=-1$ and $\beta=N+1$ we have
$$
\begin{array}{l}
-a_1^-\sum\limits_{\gamma=1}^{\infty}D(h\gamma-h)(h\gamma)+a_1^+\sum\limits_{\gamma=1}^{\infty}D(h(N+\gamma)+h)(h\gamma)=\\
=-\sum\limits_{\gamma=0}^ND(h\gamma+h)f(h\gamma,t)-\frac{1}{12}g_0\sum\limits_{\gamma=1}^{\infty}D(h\gamma-h)(h\gamma)^3-\\
-\frac{1}{4}g_1\sum\limits_{\gamma=1}^{\infty}D(h\gamma-h)(h\gamma)^2-\frac{1}{12}g_0\sum\limits_{\gamma=1}^{\infty}D(h(N+\gamma)+h)(1+h\gamma)^3+\\
+\frac{1}{4}g_1\sum\limits_{\gamma=1}^{\infty}D(h(N+\gamma)+h)(1+h\gamma)^2-\bigg( f(1,t)-\frac{1}{12}g_0+\frac{1}{4}g_1\bigg)\times\\
\times\sum\limits_{\gamma=1}^{\infty}D(h(N+\gamma)+h)-f(0,t)\sum\limits_{\gamma=1}^{\infty}D(h\gamma-h),
\end{array}\eqno(3.20)
$$
$$
\begin{array}{l}
-a_1^-\sum\limits_{\gamma=1}^{\infty}D(h(N+\gamma)+h)(h\gamma)+a_1^+\sum\limits_{\gamma=1}^{\infty}D(h\gamma-h)(h\gamma)=\\
=-\sum\limits_{\gamma=0}^ND(h(N+1)-h\gamma)f(h\gamma,t)-\frac{1}{12}g_0\sum\limits_{\gamma=1}^{\infty}D(h(N+\gamma)+h)(h\gamma)^3-\\
-\frac{1}{4}g_1\sum\limits_{\gamma=1}^{\infty}D(h(N+\gamma)+h)(h\gamma)^2-\frac{1}{12}g_0\sum\limits_{\gamma=1}^{\infty}D(h\gamma-h)(1+h\gamma)^3+\\
+\frac{1}{4}g_1\sum\limits_{\gamma=1}^{\infty}D(h\gamma-h)(1+h\gamma)^2-\bigg( f(1,t)-\frac{1}{12}g_0+\frac{1}{4}g_1\bigg)\times\\
\times\sum\limits_{\gamma=1}^{\infty}D(h\gamma-h)-f(0,t)\sum\limits_{\gamma=1}^{\infty}D(h(N+\gamma)+h).
\end{array}\eqno(3.21)
$$
Thus, for unkowns  $a_1^-,\ a_1^+,\ a_0^-, a_0^+$ we obtained the system of linear equations
(3.18)-(3.21). Solving this system we obtain  (3.14)-(3.17). This means we obtained the explicit form of the function $u(h\beta)$.

Next, using  (2.3) and (3.7), from (3.6) calculating the convolution  $C([\beta],t)=hD_2(h\beta)\ast u(h\beta)$  for $\beta=\overline{0,N}$ we respectively get
(3.11)-(3.13). Theorem is proved. \hfill $\Box$ \\

%
%\label{Akhmed2}
%
%%\endinput

\end{document}